\UseRawInputEncoding
\documentclass[compress,final,3p,times,12pt]{elsarticle}
\usepackage{amsmath}
\usepackage{amsmath,amssymb, amsthm, latexsym, mathrsfs,amsfonts, epsfig, color,graphicx,}
\usepackage{lineno,hyperref}
\renewcommand{\a}{\alpha}
\newcommand{\D}{\Delta}

\renewcommand{\l}{\left}
\renewcommand{\r}{\right}
\newcommand{\e}{\epsilon}

\newtheorem{theorem}{Theorem}
\newtheorem{lemma}{Lemma}
\newtheorem{remark}{Remark}

\modulolinenumbers[5]

\journal{Journal of Applied Mathematics Letters}

\begin{document}
\begin{frontmatter}
\title{Convergence rate for homogenization of a nonlocal model with oscillating coefficients}
  %   \author{Li Lin,  Meihua Yang\\
%\small{ School of Mathematics and Statistics, Huazhong University of Science and Technology} \\
%\small {Wuhan, 430074, China}\\
%\small{ Center for Mathematical Sciences, Huazhong University of Science and Technology} \\
%\small {Wuhan, 430074, China}\\
%Jinqiao Duan \footnote{Corresponding author: duan@iit.edu.}\\
%\small { Department of Applied Mathematics, Illinois Institute of Technology}\\
%\small {Chicago, Illnois 60616, USA}\\}
\author{Li Lin\footnote{School of Mathematics and Statistics \& Center for Mathematical Sciences, Huazhong University of Science and Technology,Wuhan, 430074, China. E-mail address: linli@hust.edu.cn   }
 \quad \&\quad Jinqiao Duan\footnote{Corresponding author. Department of Applied Mathematics, Illinois Institute of Technology, Chicago, Illnois 60616, USA. E-mail address: duan@iit.edu.}}

\begin{abstract}
This letter deals with homogenization of a nonlocal model with L\'{e}vy-type operator of rapidly oscillating coefficients. This nonlocal model describes mean residence time and other escape phenomena  for stochastic dynamical systems with non-Gaussian L\'evy noise. We derive an effective model with a specific  convergence rate. This enables efficient analysis and simulation of escape phenomena under non-Gaussian fluctuations.
\end{abstract}

\begin{keyword}
 Convergence rate, mean residence time, escape phenomena,  L\'{e}vy noise, nonlocal elliptic equations
 \MSC[2020] 34A08\sep  74Q10 \sep 41A25
 \end{keyword}
\end{frontmatter}

\section{\bf Introduction}
We  consider the homogenization of the following nonlocal partial differential equation
\begin{equation}\label{3}
\begin{cases}
  \mathcal{A}^\e u_\e(x)=f(x), \qquad x\in D,\\
  u_\epsilon|_{D^c}=g(x),
\end{cases}
\end{equation}
arising  in the study of escape phenomena of stochastic dynamical systems under L\'{e}vy fluctuations \cite{Duan}. Here     $D$ is a bounded domain   in $\mathbb{R}^d$. When $f=-1 $ and $g=0$, the solution of this equation is the mean residence time for such a stochastic system in domain $D$. Specifically,   the nonlocal operator depends on a small positive scale parameter $\varepsilon$ and is defined by
$\mathcal{A}^\e u=\frac{1}{2}\mathcal{D}(\Theta^\e\mathcal{D}^*u),$
where the coefficient  $\Theta^\e(x,z)=\Theta(\frac{x}{\e},\frac{z}{\e})$ is of period $1$ in $(x,z)$ and for a certain positive constant $\lambda$,
$\frac{1}{\lambda}<\Theta^\e(x,z)<\lambda.$

The nonlocal linear operator $\mathcal{D}$ and its adjoint operator $\mathcal{D}^*$ are defined as follows\cite{DQ}.
Given the antisymmetric kernel function
$\gamma(x,z)=(z-x)\frac{1}{|z-x|^{\frac{2+d+\a}{2}}}$ (note that $\gamma(z,x)=-\gamma(x,z)$), the nonlocal divergence $\mathcal{D}$   is defined by
$\mathcal{D}(\beta)(x):=\int_{\mathbb{R}^d}(\beta(x,z)+\beta(z,x))\cdot\gamma(x,z)dz,$  $x\in D.$
The adjoint operator $\mathcal{D^{*}}$   of  $\mathcal{D}$ is then given by
$\mathcal{D^{*}}(\phi)(x,z)=-(\phi(z)-\phi(x))\gamma(x,z)$ for $x,z\in D.$

By the way, if $\Theta \equiv 1$, then we see a relation with the nonlocal Laplace operator
$\frac{1}{2}\mathcal{D}\mathcal{D^{*}}=-(-\Delta)^{\a/2}.$  The nonlocal Laplace operator $(-\Delta)^{\a/2}$ is defined by
$$(-\Delta)^{\a/2}u(x)=\int_{\mathbb{R}^d\backslash \{x\}}\frac{u(z)-u(x)}{|z-x|^{d+\a}}dz,$$
where the integral is in the sense of Cauchy principal value, and it is the generator of a symmetric $\a$-stable L\'{e}vy motion \cite{Duant}
%\begin{remark}
%  For a function $\upsilon(x,y)$, we define
%$$(\mathcal{D}^*_x\upsilon)(x,z,y)=-(\upsilon(z,y)-\upsilon(x,y))\gamma(x,z)$$
%and
%\[
%\begin{split}
%(\mathcal{D}_x\mathcal{D}^*_x\upsilon)(x,y)&=2\int_{\mathbb{R}}
%-(\upsilon(z,y)-\upsilon(x,y))\gamma^2(x,z)dz\\
%&=-2(-\Delta)^{\a/2}_x\upsilon(x,y).
%\end{split}
%\]
%\end{remark}

\section{Main result}

\subsection{Function spaces}
In this subsection, we define the following two important spaces. One is the space $\mathcal{V}^D$, which is the counterpart of the classical Sobolev space $H^1(D)$.  The other is the space $\mathcal{X}^D$, as an analogue of the usual trace space $H^{1/2}(\partial D).$

 We set $\nu(x,y)=\gamma^2(x,z)$.  we introduce a quadratic form\cite{DS,FM,MV}  for $u: \mathbb{R}^d\rightarrow\mathbb{R}$
$$\mathcal E_D(u,u)=\frac{1}{2}\int_{\mathbb{R}^d\times\mathbb{R}^d\backslash D^c\times D^c}(u(x)-u(z))^2\nu(x,z)dxdz.$$
Now we define the corresponding Sobolev space \cite{RK}:
$$\mathcal{V}^D=\{u: \mathbb{R}^d\rightarrow\mathbb{R}\; \text{such that}\; \mathcal E_D(u,u)<\infty\}.$$
We also denote
$$\mathcal{V}^D_0=\{u\in\mathcal{V}^D; u=0\;\;\;   a.e.\; on\; D^c\}.$$
More precisely,   $\mathcal{V}^D_0=H_0^{\a/2}(D)$ from
\cite[Theorem  3.33]{MW}.  The space $H_0^{\a/2}(D)$ is the standard fractional Sobolev space.
For $s < 0$, we denote $H^s(D)$ as the dual space of $H_0^{-s}(D)$.
It is clear that   $\mathcal{V}^D $ is a Hilbert subspace of  $L^2(D)$, with the induced norm
$||u||^2_{\mathcal{V}^D}=||u||^2_{L^2(D)}+\mathcal{E}_{D}(u,u).$

We will use the notation $\mathbb{T}^d$ to denote the $d$-dimensional unit torus. The space $H^{\a/2}_\#(\mathbb{T}^d)$ of $1$-periodic functions $u\in H^{\a/2}$ such that $\int_{\mathbb{T}^d}u(y)dy=0$ will be interest in this study. Provided with the norm,
$\|u\|_{H^{\a/2}_\#(\mathbb{T}^d)}=\l(\int_{\mathbb{T}^d}\int_{\mathbb{T}^d}\frac{|u(y)-u(\eta)|^2}{|y-\eta|^{d+\a}}dyd\eta\r)^{\frac{1}{2}}.$
%$$H_\#(Y)=\{u\in H_{per}(Y)| \int_Yu(y)dy=0\}.$$

We let $$\mathcal E_D(u,v)=\frac{1}{2}\int_{\mathbb{R}^d\times\mathbb{R}^d\backslash D^c\times D^c}(u(x)-u(z))(v(x)-v(z))\nu(x,z)dxdz,$$
and
$$a^\e(u,v)=\frac{1}{2}\int_{\mathbb{R}^d\times\mathbb{R}^d\backslash D^c\times D^c}\Theta^\e(x,z)(u(x)-u(z))(v(x)-v(z))\nu(x,z)dxdz,$$
if the integrals are absolutely convergent, in particular for $u,v\in \mathcal{V}^D.$

By a solution of (\ref{3}) we mean a weak solution, which is defined as every function $u_\e\in \mathcal{V}^D$ equal to $g$ a.e. on $D^c$ such that
for every $\phi\in\mathcal{V}^D_0,$
$a^\e(u_{\e},\phi)=\int f\phi dx,$ this integral is infinite, e.g., if $D$ is bounded and $f\in L^2(D).$

\begin{remark}
  Without loss of generality, we take the coefficient $\Theta^{\e}(x,z)$ to be a symmetric function. In fact, we can define the symmetric and anti-symmetric parts of $\Theta^{\e}$:
  $$\Theta^{\e}_s(x,z)=\frac{1}{2}(\Theta^{\e}(x,z)+\Theta^{\e}(z,x)) \quad \text{and} \quad \Theta^{\e}_a(x,z)=\frac{1}{2}(\Theta^{\e}(x,z)-\Theta^{\e}(z,x)).$$
By \cite[Lemma 3.2]{RA},
 $\l(\Theta^{\e}_a(x,z)\mathcal{D}^*u,\mathcal{D}^*\phi\r)=0.$  Then $\mathcal{D}(\Theta^{\e}(x,z)\mathcal{D}^*u)=\mathcal{D}(\Theta^{\e}_s(x,z)\mathcal{D}^*u).$
\end{remark}

Let $G^\e_D(x,y)$ and $P^\e_D(x,y)$ be the Green function and Poisson kernel of $D$ for $\mathcal{A^\e}$ respectively. At the same time
we let $G_D(x,y)$ and $P_D(x,y)$ be the Green function and Poisson kernel for $\mathcal{A^\e}$ with $\Theta^\e=1$  on $D$, respectively. In this case
$\mathcal{A^\e}u=-2(-\D)^{\frac{\a}{2}}u.$

For $g: D^c\rightarrow \mathbb{R}$, we let $P^\e_D[g](x)=g(x)$ for $x\in D^c$ and
$P^\e_D[g](x)=\int_{D^c}g(y)P^\e_D(x,y)dy$ for$  x\in D.$\\
Furthermore, $u=P^\e_D[g]$ is the unique solution of the following homogeneous Dirchlet problem:
\begin{equation}\label{4}
\begin{cases}
  \mathcal{A}^\e u_\e=0, \qquad x\in D,\\
  u^\epsilon|_{D^c}=g.
  \end{cases}
\end{equation}
\begin{remark}[\cite{FM}]
  For $f\in H^{-\a/2}(D)  \text{ and }    g\in \mathcal{X}^D$,  we have the existence and uniqueness of equations (\ref{3}) and (\ref{4}) in $\mathcal{V}^D$.
\end{remark}

Next, for $\omega,\xi \in D^c,$ we let
$\gamma_D(\omega,\xi)=\int_D\int_D\nu(\omega,x)G_D(x,y)\nu(y,\xi)dxdy=\int_D\nu(\xi,x)P_D(x,\omega)dx.$\\
For $g: D^c\rightarrow \mathbb{R}$ we introduce a quadratic form
$\mathcal{H}_D(g,g)=\frac{1}{2}\int_{D^c\times D^c}(g(\omega)-g(\xi))^2\gamma_D(\xi,\omega)d\omega d\xi.$\\
Then we   define as in \cite{RK} a new space
$\mathcal{X}^D=\{g: D^c\rightarrow\mathbb{R} \; \text{such that}\; \mathcal{H}_D(g,g)<\infty\}.$

If  $g\in\mathcal{X}^D$ and $x\in D,$ we   obtain $\int_{D^c}g(z)^2P_D(x,z)dz<\infty.$ We fix an  arbitrary (reference) point $x_0\in D.$ For $g\in \mathcal{X}^D,$ we let
$|g|^2_{D^c}=\int_{D^c}g(z)^2P_D(x_0,z)dz$
(we omit $x_0$ from the notation).
Then  $\mathcal{X}^D$ is a Hilbert space with the induced norm
$||g||^2_{\mathcal{X}^D}=|g|^2_{D^c}+\mathcal{H}_D(g,g).$

%For every positive real number $m$,   write $m=\sigma+\tau,$ with $\sigma$
%  an integer and $\tau\in (0,1).$ We define a Sobolev space
%$$H^m(D)=\{u\in H^m: D^\sigma u\in H^{\tau}(D)\}.$$

\subsection{Effective equation and convergence rate}

Our main result is that the heterogeneous model (\ref{3}) is approximated by a homogenized effective model (\ref{2}) below, with convergence rate $\frac12$. This is stated in the following theorem.

\begin{theorem}\label{T1}
  For scale parameter $\varepsilon$ sufficiently small, the solution $u_\e$ of heterogeneous model (\ref{3}) is approximated   by the solution $u_0$ of the following homogenized equation
  \begin{align}\label{2}
  \begin{cases}
    -a_1(-\Delta)^{\a/2}u_0-a_2\mathcal{F}u_0(x) =f(x), \qquad x\in D,\\
    u_0|_{D^c}=g(x).
  \end{cases}
\end{align}
with coefficients
\begin{align*}
  &a_1=\int_{\mathbb{T}^d\times\mathbb{T}^d}\Theta(y,\eta)dyd\eta,\qquad
  a_2=\frac{1}{2}\int_{\mathbb{T}^d\times\mathbb{T}^d}\Theta(y,\eta)\mathcal{D}_y^*\chi(y) dyd\eta,\\
  &\zeta(u_0)(x)=\frac{1}{|D|}\int_{D}(D^*u_0)(x,z)dz,\qquad
  \mathcal{F}u_0(x)=\mathcal{D}|_{D}\zeta(u_0)(x) = \int_{D}\l[\zeta(u_0)(x)+\zeta(u_0)(z)\r]\gamma(x,z)dz,
\end{align*}
where the function $\chi(y)$   is the unique solution of the following variational problem
\begin{equation}\label{ka1}
\begin{cases}
 $${\hat{a}}(\chi,v)=\int_{\mathbb{T}^d\times\mathbb{T}^d}\Theta(y,\eta)D^*_yvdy d\eta,\\
 \chi\in H_{\#}^{\a/2}(\mathbb{T}^d).\\
\end{cases}
\end{equation}
Moreover, if $f\in C^\infty(\overline{D})$ and $\chi\in L^\infty(\mathbb{T}^d)$, then $u_\e$ has the following asymptotic expansion
$$u_{\e}=u_0-\varepsilon^{\frac{1+\a}{2}}\frac{1}{|D|}\int_{D}(D^*u_0)(x,z)dz\cdot\chi(\frac{x}{\e})+ R(\varepsilon),$$
and there exists
a constant $C$ (independent of $\varepsilon$) such that the remainder is estimated as
$$\l\|R(\varepsilon)\r\|_{\mathcal{V}^D}=\l\|u_{\e}-(u_0-\varepsilon^{\frac{1+\a}{2}}\frac{1}{|D|}\int_{D}(D^*u_0)(x,z)dz\cdot\chi(\frac{x}{\e}))\r\|_{\mathcal{V}^D}\leq C\varepsilon^{1/2}.$$
This says that $u_{\e}\rightarrow u_0,$ in Sobolev space $\mathcal{V}^D,$ with convergence rate $\frac{1}{2}.$
\end{theorem}

In order to prove this theorem, we recall some lemmas. The next result is due to   \cite{RK}.
\begin{lemma}\label{l1}
  Let $D\subset \mathbb{R}^d$ be bounded, open and Lipschitz, $|\partial D|=0.$
  \begin{itemize}
    \item  If $g\in \mathcal{X}^D,$ then $P^\e_D[g]\in \mathcal{V}^D$ and $\mathcal{E}_D(P^\e_D[g],P^\e_D[g])=\mathcal{H}_D(g,g).$
    \item  If $u\in \mathcal{V}^D, $ then $g=u|_{D^c}\in \mathcal{X}^D$ and $\mathcal{E}_D(u,u)\geq\mathcal{H}_D(g,g).$
  \end{itemize}
\end{lemma}

Let $X=\{X_t\}_{t\geq0}$ be a L\'{e}vy process with $(0,\nu,0)$ as the L\'{e}vy triplet on a probability space $(\Omega,\mathcal{F},\mathbb{P}).$ We introduce the time of the first exit of $X$ from $D,$
$\tau_D=\tau_D(X)=\inf\{t\geq0: X_t\notin D\}.$

\begin{lemma}\label{l2}
  The assumptions are the same as in Lemma \ref{l1}. Then for every $g\in \mathcal{X}^D,$ there exists a positive constant $C(D)$ such that $||P^\e_D[g]||_{\mathcal{V}^D}\leq C(D)||g||_{\mathcal{X}^D}.$
  \begin{proof}
    We write $U\subset\subset D$ if U is an open set, its closure $\overline{U}$
    is bounded, and $\overline{U}\subset D.$ Let $\tilde{u}^\e$ be the unique solution of the homogeneous Dirchlet problem (\ref{4}) and $\mathbb{E}^x$ be the expectation for
     $X_t$ start at $x\in D$. We have
     $\mathbb{E}^x \tilde{u}_\e(X_{\tau_D})=\int_{D^c}\tilde{u}_\e(y)P_D^\e(x,y)dy.$
     Then we obtain $$||P^\e_D[g]||^2_{L^2(D)}=\int_{D}(\int_{D^c}g(y)P_D^\e(x,y)dy)^2dx\leq\int_{D}\int_{D^c}g(y)^2P_D^\e(x,y)dydx,$$
     due to the fact that $\int_{D^c}P_D^\e(x,y)dy=1,$ for $x\in D.$
     That is to say $||P^\e_D[g]||^2_{L^2(D)}\leq\mathbb{E}^x \tilde{u}^2_\e(X_{\tau_D}).$

     Note that $\mathbb{E}^x\tilde{u}_\e(X_{\tau_D})$ is a closed martingale\cite[Remark 4.4]{RK} and for $x\in U,$ $x\rightarrow\int_{U^c}g(y)^2P_U^\e(x,y)dy$ satisfies the Harnack inequality\cite{BK}. Thus
     \begin{equation*}
     \begin{split}
       ||P^\e_D[g]||^2_{L^2(D)}&\leq \lim_{x\in U\subset\subset D}\int_U\int_{U^c}g(y)^2P_U^\e(x,y)dydx
       \leq\lim_{x\in U\subset\subset D}\int_U
       C \int_{U^c}g(y)^2P_U^\e(x_0,y)dydx
       =C(D)|g|_{D^c}^2.
     \end{split}
    \end{equation*}
    By   Lemma \ref{l1}, we obtain
    \begin{equation*}
     \begin{split}
     ||P^\e_D[g]||^2_{\mathcal{V}^D}&=||P^\e_D[g]||^2_{L^2(D)}+\mathcal{E}_{D}(P^\e_D[g],P^\e_D[g])\leq C(D)|g|_{D^c}^2+\mathcal{H}_D(g,g)
     \leq C(D)||g||^2_{\mathcal{X}^D}.
     \end{split}
    \end{equation*}
    Hence   Lemma 2 follows.
  \end{proof}
\end{lemma}
Next, we obtain a uniform estimate concerning the solution $u_\e$ for the original heterogeneous equation (\ref{3}).
\begin{lemma}\label{le}
  Let  $f$ in $H^{-\a/2}(D)$, $g$ in $\mathcal{X}^D$ and $u_\e$ be the unique solution of the original heterogenous equation (\ref{3}).
  Then there exist two positive constants $C_1, C_2$  such that
  $$||u_\e||_{\mathcal{V}^D}\leq C_1||f||_{H^{-\a/2}(D)}+C_2||g||_{\mathcal{X}^D}.$$
  \begin{proof}
   % For $f: D\rightarrow \mathbb{R}$ we let $G^\e_D[f](x)=0\; \text{for}\; x\in D$ and
%    $$G^\e_D[f](x)=\int_Df(y)G^\e_D(x,y)dy \quad \text{for}\quad x\in D^c.$$
%    Due to Grzywny \cite{GT}, The Dirchlet problem (\ref{3}) has a unique solution:
%    $$u_\e=G^\e_D[f](x)+P^\e_D[g](x).$$
    From Lemma \ref{l2}, we can obtain the following conclusion. For every $g\in \mathcal{X}^D,$ there exists $G\in\mathcal{V}^D$ and a linear operator $\rho$ such that $\rho(G)=G|_{D^c}=g$ and
    $||G||_{\mathcal{V}^D}\leq C(D)||g||_{\mathcal{X}^D}.$
    For every $v\in\mathcal{V}_0^D,$ we have
    \begin{small}
    \begin{equation*}
     \begin{split}
     |(\mathcal{A}^\e G,v)|&=\frac{1}{2}|(\mathcal{D}(\Theta^\e\mathcal{D}^*G),v)|
     \leq \lambda C||\mathcal{D}^*G||_{L^2(D\times\mathbb{R}^d)}||\mathcal{D}^*v||_{L^2(D\times\mathbb{R}^d)}
     \leq 4\lambda C||G||_{\mathcal{V}^D}||v||_{\mathcal{V}_0^D}.
    \end{split}
    \end{equation*}
    \end{small}
    Then we obtain $\mathcal{D}(\Theta^\e\mathcal{D}^*G)\in (\mathcal{V}_0^D)^*.$
    That is to say, for $x\in D,$ we have $f-\mathcal{D}(\Theta^\e\mathcal{D}^*G)\in H^{-\a/2}.$
    Recall that
    $$a^\e(u,v)=\frac{1}{2}\int_{\mathbb{R}^d\times\mathbb{R}^d\backslash D^c\times D^c}\Theta^\e(x,z)(u(x)-u(z))(v(x)-v(z))\nu(x,z)dxdz,$$
  for every $v\in\mathcal{V}_0^D.$ We can find a unique $v_\e\in\mathcal{V}_0^D$ such that $a^\e(v_\e,v)=<f-\frac{1}{2}\mathcal{D}(\Theta^\e\mathcal{D}^*G),v>_{(\mathcal{V}_0^D)^*,\mathcal{V}_0^D}.$
  Due to the Poincar\'{e}    inequality \cite{FM}, there exists a constant $C\geq1,$ for every $u\in \mathcal{V}_0^D,$
  $$||u||_{L^2(D)}^2\leq C\int_{\mathbb{R}^d\times \mathbb{R}^d\backslash D^c\times D^c}(u(x)-u(z))^2\nu(x,z)dxdz.$$
  We thus obtain $||u||_{\mathcal{V}_0^D}\leq(2C+1)||\mathcal{D}^*u||_{L^2(D\times \mathbb{R}^d)}.$
  In other words, the space $\mathcal{V}_0^D$ can be equipped by the norm $||\mathcal{D}^*u||_{L^2(D\times \mathbb{R}^d)}.$
  Then $a^\e(v,v)\geq C||v||^2_{\mathcal{V}_0^D}$ and $|a^\e(u,v)|\leq C||u||_{\mathcal{V}_0^D}||v||_{\mathcal{V}_0^D}.$
  From the Lax-Milgram theorem,
  \begin{equation*}
    \begin{split}
      \frac{1}{\lambda}||v_\e||_{\mathcal{V}_0^D}^2\leq a^\e(v_\e,v_\e)=<f-\frac{1}{2}\mathcal{D}(\Theta^\e\mathcal{D}^*G),v_\e>_{(\mathcal{V}_0^D)^*,\mathcal{V}_0^D}
      \leq||f-\frac{1}{2}\mathcal{D}(\Theta^\e\mathcal{D}^*G)||_{H^{-\a/2}(D)}||v_\e||_{\mathcal{V}_0^D}.
      \end{split}
  \end{equation*}
  We can see that $$||v_\e||_{\mathcal{V}_0^D}\leq \lambda ||f-\frac{1}{2}\mathcal{D}(\Theta^\e\mathcal{D}^*G)||_{H^{-\a/2}(D)}.$$
  Set $u_\e=v_\e+G.$ By the linearity of $\rho$, we have $\rho(u_\e)=\rho(G)=g.$ Furthermore,
  $$a^\e(u_\e,v)=a^\e(v_\e,v)+a^\e(G,v)=(f,v),$$
   which means that $u_\e$ is the unique solution of the  original heterogenous equation (\ref{3}).
   Then
   \begin{equation*}
    \begin{split}
     ||u_\e||_{\mathcal{V}^D}&\leq||u_\e-G||_{\mathcal{V}^D}+||G||_{\mathcal{V}^D}
     \leq C||v_\e||_{\mathcal{V}_0^D}+\sqrt{C(D)}||g||_{\mathcal{X}^D}\\
     &\leq \lambda  C||f-\frac{1}{2}\mathcal{D}(\Theta^\e\mathcal{D}^*G)||_{H^{-\a/2}(D)}+\sqrt{C(D)}||g||_{\mathcal{X}^D}.
      \end{split}
  \end{equation*}
  On the other hand,
  \begin{equation*}
    \begin{split}
    <\mathcal{D}(\Theta^\e\mathcal{D}^*G), v>_{(\mathcal{V}_0^D)^*,\mathcal{V}_0^D}&=\int_{\mathbb{R}^d\times \mathbb{R}^d\backslash D^c\times D^c}
    \Theta^\e(x,z)\mathcal{D}^*G\mathcal{D}^*vdxdz
    \leq\lambda C||\mathcal{D}^*G||_{L^2(D\times\mathbb{R}^d)}||\mathcal{D}^*v||_{L^2(D\times\mathbb{R}^d)}\\
    &\leq C||G||_{\mathcal{V}^D}||\mathcal{D}^*v||_{L^2(D\times\mathbb{R}^d)}
    \leq C||g||_{\mathcal{X}^D}||v||_{\mathcal{V}_0^D}.
     \end{split}
  \end{equation*}
  That implies $||\mathcal{D}(\Theta^\e\mathcal{D}^*G)||\leq C||g||_{\mathcal{X}^D}$.
  Hence
    $||u_\e||_{\mathcal{V}^D}\leq C_1||f||_{H^{-\a/2}(D)}+C_2||g||_{\mathcal{X}^D}.$
    This completes the proof.  \end{proof}
\end{lemma}

\subsection{Proof of Theorem 1}
We are now ready to prove our main result in Theorem \ref{T1}.

\textbf{Step 1:} First, we will derive the homogenized equation for $\varepsilon$ sufficiently small.
  For a function $\upsilon(x,y)$, we define
$$(\mathcal{D}^*_x\upsilon)(x,z,y)=-(\upsilon(z,y)-\upsilon(x,y))\gamma(x,z)$$
and
\[
\begin{split}
(\mathcal{D}_x\mathcal{D}^*_x\upsilon)(x,y)&=2\int_{\mathbb{R}^d}
-(\upsilon(z,y)-\upsilon(x,y))\gamma^2(x,z)dz\\
&=-2(-\Delta)^{\a/2}_x\upsilon(x,y).
\end{split}
\]
Denote $\eta=\frac{z}{\epsilon}$ a variable on the period: $\eta \in \mathbb{T}^d.$ We look for a  formal asymptotic expansion:
$$u_\e=u_0(x,\frac{x}{\e})+\varepsilon^{\frac{1+\a}{2}}u_1(x,\frac{x}{\e})+o(\varepsilon^{\frac{1+\a}{2}}),$$
with $u_i(x,y)$, for $i=1,2$,  such that
\begin{equation*}
\begin{cases}
u_i(x,y) \quad \text{is defined for} \quad x\in{D} \quad\text{and}\quad y\in{\mathbb{T}^d},\\
u_i(\cdot,y)\quad \text{is $1$-periodic}.\\
\end{cases}
\end{equation*}
For every function $h,$ we denote $h^\e(x)=h(\frac{x}{\epsilon})$. Thus
$$\mathcal{D}^*u_\e=\mathcal{A}_1u_0+\varepsilon^{-\frac{1+\a}{2}}(\mathcal{A}_0u_0)^\e+\varepsilon^{\frac{1+\a}{2}}\mathcal{A}_1u_1+(\mathcal{A}_0u_1)
^\e+o(\varepsilon^{\frac{1+\a}{2}}),$$
where $$\mathcal{A}_0v(x,y):=(\mathcal{D}^*_yv)(x,y,\eta), \mathcal{A}_1v(x,y):=(\mathcal{D}^*_xv)(x,z,y).$$
Then we have $\mathcal{A}_0u_0=0$. Furthermore, $u_0(x,y)=u_0(x).$

Let $C_{per}(\mathbb{T}^d))$ be the subspace of $C(\mathbb{R}^n)$ of $1$-periodic functions. For every $v\in \mathcal{M}(D, C_{per}(\mathbb{T}^d))$\\($\mathcal{M}(D)$ is the space of functions in $C^\infty$ with compact support), we  denote $v_\e=v(x,\frac{x}{\e}).$  We   conclude
%$$(\mathcal{A}^\e u_\e,v_\e)=\frac{1}{2}(\Theta^\e \mathcal{D}^*u_\e, D^*v_\e)=(f,v_\e).$$ Therefore
%$$\frac{1}{2}(\Theta^\e(\mathcal{A}_1u_0+(\mathcal{A}_0u_1)^\e),D^*v_\e)=(f,v_\e).$$
\begin{equation*}
  \begin{split}
(f,v_\e)&=\frac{1}{2}(\Theta^\e(\mathcal{A}_1u_0+(\mathcal{A}_0u_1)^\e),D^*v_\e)
=\frac{1}{2}(\Theta^\e(\mathcal{A}_1u_0+(\mathcal{A}_0u_1)^\e),D^*v_\e)_{L^2(D\times D)}\\&+\int_{D\times D^c}\Theta^\e(u_0(x)-g(z))v_{\e}(x)\gamma^2(x,z)dzdx
:=I_1^\e+I_2^\e.
  \end{split}
\end{equation*}
By \cite[Lemma 2.34]{GP}, for every $v\in \mathcal{M}(D, C_{per}(\mathbb{T}^d)),$ $I_1^\e, I_2^\e$ converges, as $\varepsilon$ goes to $0.$
Now take $v_\e(x)=v_0(x)+\varepsilon^{\frac{1+\a}{2}}v_1(x,\frac{x}{\e}), $ where $v_0(x)\in L^2(D)$ and $v_1\in L^2(D, C_{per}(\mathbb{T}^d)).$
As $\varepsilon$ goes to $0,$ we have
\begin{equation}\label{lab}
I_1^\e\rightarrow\frac{1}{2}\int_{\mathbb{T}^d\times \mathbb{T}^d}\int_{D\times D}\Theta(y,\eta)(\mathcal{D}^*u_0+
\mathcal{D}_y^*u_1)(\mathcal{D}^*v_0+\mathcal{D}_y^*v_1)dxdzdyd\eta.
\end{equation}
On one hand, let $v_0=0,$ we have $I_2^\e\rightarrow 0$ and
\begin{small}
\begin{equation}\label{1}
\begin{split}
\int_{D\times D}(\Theta(y,\eta)\mathcal{D}_y^*u_1,\mathcal{D}_y^*v_1)_{L^2({\mathbb{T}^d\times\mathbb{T}^d})}dxdz&=-\int_{D\times D}(\Theta(y,\eta)\mathcal{D}^*u_0,\mathcal{D}_y^*v_1)_{L^2{(\mathbb{T}^d\times\mathbb{T}^d})}dxdz.
\end{split}
\end{equation}
\end{small}
%Let us introduce two quadratic forms $${\hat{a}}(w,v)=\frac{1}{2}\int_{\mathbb{T}^d\times\mathbb{T}^d}\Theta(y,\eta)D^*_ywD^*_yvdy d\eta$$
%for all $w,v\in H_{\#}^{\a/2}(\mathbb{T}^d)$, and
%$$a^\e(u,\phi)=\frac{1}{2}\int_{D}\int_{D}\Theta^\e\mathcal{D}^*u(x,z)\mathcal{D}^*\phi(x,z)dxdz $$
%for all $u,\phi\in H_{\#}^{\a/2}(D).$
For all $w,v\in H_{\#}^{\a/2}(\mathbb{T}^d)$, %and $u,\phi\in H_{\#}^{\a/2}(D),$
we introduce two quadratic forms:
${\hat{a}}(w,v)=\frac{1}{2}\int_{\mathbb{T}^d\times\mathbb{T}^d}\Theta(y,\eta)D^*_ywD^*_yvdy d\eta.$
%a^\e(u,\phi)=\frac{1}{2}\int_{D}\int_{D}\Theta^\e\mathcal{D}^*u(x,z)\mathcal{D}^*\phi(x,z)dxdz. $$
From equation (\ref{1}), we have
$|D|\int_{D}{\hat{a}}(u_1,v)dx=-\int_D(\int_{D}\mathcal{D}^*u_0(x,z)dz){\hat{a}}(\chi,v)dx,$
where $\chi(y)$ is the unique solution of the following variational problem
\begin{equation}\label{ka}
\begin{cases}
 {\hat{a}}(\chi,v)=\int_{\mathbb{T}^d\times\mathbb{T}^d}\Theta(y,\eta)D^*_yvdy d\eta,\\
 \chi\in H_{\#}^{\a/2}(\mathbb{T}^d),\\
\end{cases}
\end{equation}
for all $v\in H_{\#}^{\a/2}(\mathbb{T}^d).$
So
  \begin{equation}\label{u1}
  u_1(x,y)=-\frac{1}{|D|}\int_{D}(D^*u_0)(x,z)dz\cdot\chi(y)\in L^2(\mathbb{R}^d,H_{\#}^{\a/2}(\mathbb{T}^d)).
  \end{equation}

Moreover, let $v_1=0,$ we conclude that
$$I_2\rightarrow \frac{1}{2}\int_{\mathbb{T}^d\times\mathbb{T}^d}\Theta(y,\eta)dyd\eta\int_{D\times D^c}(u_0(x)-g(z))v_0(x)\nu(x,z)dxdz,$$
as $\varepsilon$ goes to $0.$
Substituting the representation of $u_1$ in (\ref{u1}) into the equation (\ref{lab}), we have
 \begin{align}
  \begin{cases}
    -a_1(-\Delta)^{\a/2}u_0-a_2\mathcal{F}u_0(x)=f(x), \qquad x\in D,\\
    u_0|_{D^c}=g(x).
  \end{cases}
\end{align}
Here
\begin{align*}
  &a_1=\int_{\mathbb{T}^d\times\mathbb{T}^d}\Theta(y,\eta)dydn,\qquad
  a_2=\frac{1}{2}\int_{\mathbb{T}^d\times\mathbb{T}^d}\Theta(y,\eta)\mathcal{D}_y^*\chi dydn,\\
  &\zeta(u_0)(x)=\frac{1}{|D|}\int_{D}(D^*u_0)(x,z)dz,\quad
  \mathcal{F}u_0(x)=\mathcal{D}|_{D}\zeta(u_0)(x) = \int_{D}\l[\zeta(u_0)(x)+\zeta(u_0)(z)\r]\gamma(x,z)dz.
\end{align*}

%
%  \begin{equation}
%    ||u^{\e}-(u_0-\varepsilon^{\frac{1+\a}{2}}\frac{1}{|D|}\int_{D}(D^*u_0)(x,z)dz\cdot\chi(\frac{x}{\e}))||_{\mathcal{V}^D}\leq C_4\varepsilon^{(\a-1)/2}
%  \end{equation}
\textbf{Step 2:} In this step, we use the letter
$C$ for a constant independent of $\varepsilon$. We will prove the convergence rate to be $\frac12$, in the Sobolev space $\mathcal{V}^D.$

    Setting
    $$Z_\e(x)=u_{\e}(x)-(u_0+\varepsilon^{\frac{1+\a}{2}}u_1)(x,\frac{x}{\e}),$$
    we have
\begin{equation}
  \begin{cases}
   \mathcal{A}_{\e}Z_\e=-\frac{1}{2}\varepsilon^{\frac{1+\a}{2}}\mathcal{D}(\Theta\mathcal{D}_x^*u_1)^\e:=\frac{1}{2}\varepsilon^{\frac{1+\a}{2}}F_\e(x)\quad\text{in} \; D,\\
   Z_\e|_{D^c}=-\varepsilon^{\frac{1+\a}{2}}u_1(x,\frac{x}{\epsilon}):=\varepsilon^{\frac{1+\a}{2}}K_\e(x)
   :=\varepsilon^{\frac{1+\a}{2}}K(x,\frac{x}{\epsilon})
    \quad\text{on} \; D^c.
  \end{cases}
\end{equation}
We can easily check that $F_\e\in H^{-\a/2}(D).$
Let us now look at the function $K_\e.$  We prove the following estimate:
$||K_\e||_{\mathcal{X}^D}\leq C \varepsilon^{-\a/2}.$
%From the definition of $u_1,$  if $|x|\rightarrow \infty,$ we have $K_\e(x)\sim \frac{1}{|x|^{\frac{1+\a}{2}}}.$
%Then there exists a positive constant $N$ such that  $\int_{\{|x|>N,|z|>N\}}(K_\e(x)-K_\e(z))^2
%\nu(x,z)dxdz=o(\epsilon).$

For a large enough constant $N,$ we set $M^c=\{x, |x|>N\}.$ Introduce the function $m_\e$ defined as follows:
%\begin{equation*}
%  \begin{cases}
%   m_\e=1 \quad \text{if} \quad dist(x,\partial D)\leq\varepsilon, x\in D \; \text{or} \; x\in D^c/M^c \; \text{or} \; dist(x,\partial M)\leq\varepsilon, x\in M^c\\
%   m_\e=0 \quad \text{if} \quad dist(x,\partial D)\geq 2\varepsilon, x\in D\;\text{or} \; dist(x,\partial M)\geq 2\varepsilon, x\in M^c\\
%   ||\nabla m_\e||_{L^\infty(\mathbb{R}^d)}\leq \epsilon^{-1}C,
%  \end{cases}
%\end{equation*}

\begin{equation*}
  \begin{cases}
   m_\e=1 \quad \text{if} \quad dist(x,\partial D)\leq\varepsilon, x\in D \; \text{or} \; x\in D^c/M^c \\
   m_\e=0 \quad \text{if} \quad dist(x,\partial D)\geq 2\varepsilon, x\in D\;\text{or} \; x\in M^c\\
   ||\nabla m_\e||_{L^\infty(\mathbb{R}^d)}\leq \varepsilon^{-1}C,
  \end{cases}
\end{equation*}

Moreover, we set $m^{\e}\in C^\infty(\mathbb{R}^d/\partial M),$ and the derivative of the function $m^{\e}$ at $\{x\in D:dist (x,\partial D)=2\varepsilon\}$ is $0.$
Set $\psi_\e=m_\e K_\e.$ The support of $\psi_\e$ in the domain $D$ is a neighbourhood of thickness $2\varepsilon$ which we denote by $U^\e.$
% Due to the definition of $m_\e,$ the support of $\psi_\e$

First of all, we prove the estimates $\l\|\psi_\e\r\|_{\mathcal{V}^{U^\e}}\leq C\varepsilon^{-\a/2}.$

Clearly, from the definition of $m_\e$ and the regularity properties of $u_0,$ we have
$\l\|\psi_\e\r\|_{L^2(U^\e)}\leq C.$
Moreover,
we have
\begin{small}
\begin{equation*}
\begin{split}
\mathcal{D}^*\psi_\e&=\mathcal{D}^*(m_\e K_\e)=((m_\e K_\e)(x)-(m_\e K_\e)(z))\gamma(x,z)\\
&=m_\e(x)\l((\mathcal{D}^*_xK)(x,z,\frac{x}{\epsilon})+\e^{-\frac{1+\a}{2}}(\mathcal{D}^*_yK)^\e|_{x=z}(z,\frac{x}{\epsilon},\frac{z}{\epsilon})\r)+
(\mathcal{D}^*m_\e)(x,z)\cdot K_\e(z).
\end{split}
\end{equation*}
\end{small}
That is to say \cite[Lemma 3.2]{RA},
\begin{small}
\begin{equation*}
\begin{split}
\mathcal E_{U^\e}(\psi_\e,\psi_\e)&\leq C\int_{U^\e\times \mathbb{R}^d }(\mathcal{D}^*\psi_\e)^2dzdx
\leq \int_{U^\e\times \mathbb{R}^d}(m_\e)^2(x)(\mathcal{D}^*_xK)^2(x,z,\frac{x}{\epsilon})dzdx\\
&+\varepsilon^{-1-\a}\int_{U^\e\times \mathbb{R}^d}(m_\e)^2(x)
\l((\mathcal{D}^*_yK|_{x=z})^\e(z,\frac{x}{\epsilon},\frac{z}{\epsilon})\r)^2dzdx\\
&+\int_{U^\e\times \mathbb{R}^d}(\mathcal{D}^*m_\e)^2(x,z)\cdot (K_\e)^2(z)dxdz
:=J_1+J_2+J_3.
\end{split}
\end{equation*}
\end{small}
We set $K(x,\frac{x}{\e})=h(x)\chi(\frac{x}{\e}),$ here $h(x)=\frac{1}{|D|}\int_{D}(D^*u_0)(x,z)dz.$
Clearly, from the regularity of properties of  $u_0$, $\chi$, we have
\begin{align*}
J_1&=\int_{U^\e\times \mathbb{R}^d}(m_\e)^2(x)(\mathcal{D}^*_xh)^2(x,z)\chi^2(\frac{x}{\e})dzdx
\leq C\int_{U^\e\times \mathbb{R}^d}(\mathcal{D}^*_xh)^2(x,z)dzdx\leq C,\\
J_2&=\varepsilon^{-1-\a}\int_{U^\e\times \mathbb{R}^d}(m_\e)^2(x)\l((\mathcal{D}^*_y\chi)^\e\r)^2h^2(z)dzdx
=\int_{U^\e\times \{\mathbb{R}^d\cap|x-z|\geq1\}}(m_\e)^2(x)\frac{(\chi^{\e}(x)-\chi^{\e}(z))^2}{|x-z|^{d+\a}}h^2(z)dzdx\\
&+\varepsilon^{-1-\a}\int_{U^\e\times \{\mathbb{R}^d\cap|x-z|\leq1\}}(m_\e)^2(x)\l((\mathcal{D}^*_y\chi)^\e\r)^2h^2(z)dzdx
\leq C\int_{U^\e\times \{\mathbb{R}^d\cap|x-z|\geq1\}}(m_\e)^2(x)h^2(z)dzdx\\
&+\varepsilon^{-1-\a}\int_{U^\e\times \{\mathbb{R}^d\cap|x-z|\leq1\}}(m_\e)^2(x)\l((\mathcal{D}^*_y\chi)^\e\r)^2h^2(z)dzdx
%&\leq C+C \e^{-1-\a}\int_{U^\e\times \{\mathbb{R}^d\cap|x-z|\leq1\}}h^2(z)(m_\e)^2(x)\int_{\mathbb{T}^d\times \mathbb{T}^d}(\mathcal{D}^*_y\chi)^2(y,\eta)dyd\eta dzdx\\
\leq C \varepsilon^{-1-\a}\int_{U^\e\times D}(\mathcal{D}^*u_0)^2dzdx,\\
J_3&=\int_{U^\e\times \mathbb{R}^d}(\mathcal{D}^*m_\e)^2(x,z)\cdot (K_\e)^2(z)dzdx
=\int_{U^\e\times \{\mathbb{R}^d\cap|x-z|\leq1\}}\frac{[(m^\e)'(\xi)]^2(x-z)^2}{|x-z|^{d+\a}}h^2(z)\chi^2(\frac{z}{\epsilon})dzdx\\
&+C\int_{U^\e\times \{\mathbb{R}^d\cap|x-z|\geq1\}}\frac{1}{|x-z|^{d+\a}}h^2(z)\chi^2(\frac{z}{\epsilon})dzdx
\leq C \varepsilon^{-2}\int_{U^\e\times D}(\mathcal{D}^*u_0)^2dzdx.
\end{align*}
%Since the fact $\chi\in C^1(\mathbb{T}^d),$ we have \cite[Theorem 2.28]{GP}
%\begin{small}
%$$\int_{\mathbb{R}\times U^\e}(\mathcal{D}^*m_\e)(x,z)\cdot K(z,\frac{z}{\epsilon})dzdx\leq\l(\int_{\mathbb{R}\times U^\e}
%(\mathcal{D}^*m_\e)^2(x,z)dzdx\r)^{1/2}\l(\int_{U^\e}K^2(z,\frac{z}{\epsilon})dz\r)^{1/2}.$$
%\end{small}
That is to say
%We thus conclude that
$\mathcal E_{U^\e}(\psi_\e,\psi_\e)\leq J_1+J_2+J_3\leq C\varepsilon^{-1-\a}\int_{U^\e\times D}(\mathcal{D}^*u_0)^2dzdx.$\\
We can  use   a result from \cite{Oleinik}(Chapter $1$, Lemma $1.5$), which states that there exists  positive constants $C$, independent of $\varepsilon,$ such that
$\l\|u_0\r\|_{L^2(U^\e)}\leq C\varepsilon^{\frac{1}{2}}\l\|u_0\r\|_{H^1(D)}.$ We can conclude
$$\int_{U^\e}\l(\int_{D}(\mathcal{D}^*u_0)^2dx\r)dz\leq C\varepsilon\l\|u_0\r\|^2_{H^1(D)}.$$
Then $\l\|\psi_\e\r\|^2_{\mathcal{V}^{U^\e}}=\l\|\psi_\e\r\|^2_{L^2(U^\e)}+\mathcal E_{U^\e}(\psi_\e,\psi_\e)\leq C\varepsilon^{-\a}.$
That is to say
\begin{equation}\label{yao3}
\l\|\psi_\e\r\|_{\mathcal{V}^{U^\e}}\leq C\varepsilon^{-\a/2}.
\end{equation}

Secondly, we will show $\l\|K_\e\r\|_{\mathcal{X}^D}=\l\|\psi_\e\r\|_{\mathcal{X}^D}+C.$
Observe now that $\psi_\e=K_\e$ on $D^c/M^c$ and $\psi_\e=0$ on $M^c.$ Then
\begin{align*}
\l\|K_\e\r\|^2_{\mathcal{X}^D}&=C \int_{D^c/M^c\times M^c}(K_\e(x)-K_\e(z))^2\gamma_D(x,z)dzdx
+\int_{D^c/M^c\times D^c/M^c}(K_\e(x)-K_\e(z))^2\gamma_D(x,z)dzdx\\
&+\int_{M^c\times M^c}(K_\e(x)-K_\e(z))^2\gamma_D(x,z)dzdx,\\
\l\|\psi_\e\r\|^2_{\mathcal{X}^D}&=C \int_{D^c/M^c\times M^c}(\psi_\e(x)-\psi_\e(z))^2\gamma_D(x,z)dzdx
+\int_{D^c/M^c\times D^c/M^c}(K_\e(x)-K_\e(z))^2\gamma_D(x,z)dzdx.
\end{align*}
Recall that \cite[Theorem 2.6]{RK}
\begin{equation*}
\gamma_D(x,z)\approx\left\{
\begin{aligned}
&\nu(\delta_D(z))\nu(\delta_D(x))\quad \text{if} \quad diam(D)\leq\delta_D(x),\delta_D(z),\\
&\nu(\delta_D(z))/V(\delta_D(x))\quad \text{if} \quad \delta_D(x)\leq diam(D)\leq\delta_D(z),\\
&\frac{\nu(r(x,z))V^2(r(x,z))}{V(\delta_D(x))V(\delta_D(z))} \quad \text{if}\quad \delta_D(x),\delta_D(z)\leq diam(D),
\end{aligned}
\right.
\end{equation*}
here $\delta_D(x)=dist(x,\partial D), r(x,y)=\delta_D(x)+|x-y|+\delta_D(x)$ and $V(r)=C r^{\a/2}.$
Then, we have
\begin{small}
\begin{equation*}
\begin{split}
&\int_{M^c\times M^c}(K_\e(x)-K_\e(z))^2\gamma_D(x,z)dzdx
\leq C\int_{M^c\times M^c}(K_\e(x)-K_\e(z))^2\l(\delta_D(x)\r)^{-1-\a}\l(\delta_D(z)\r)^{-1-\a}dzdx
\leq C,
\end{split}
\end{equation*}
\end{small}
 and
\begin{small}
\begin{equation*}
\begin{split}
&\int_{D^c/M^c}\int_{\times M^c}(K_\e(x)-K_\e(z))^2\gamma_D(x,z)dzdx\\
&=C \int_{D^c/M^c\cap\{x:\delta_D(x)\leq diam(D)\}}\int_{ M^c}(K_\e(x)-K_\e(z))^2\l(\delta_D(x)\r)^{-\a/2}\l(\delta_D(z)\r)^{-1-\a}dzdx\\
&+C \int_{D^c/M^c\cap\{x:\delta_D(x)\geq diam(D)\}}\int_{M^c}(K_\e(x)-K_\e(z))^2\l(\delta_D(x)\r)^{-1-\a}\l(\delta_D(z)\r)^{-1-\a}dzdx\\
&\leq C.
\end{split}
\end{equation*}
\end{small}
Then, we get the conclusion
\begin{equation}\label{yao1}
\l\|K_\e\r\|^2_{\mathcal{X}^D}=\l\|\psi_\e\r\|^2_{\mathcal{X}^D}+C.
\end{equation}

Next, we will show $\l\|\psi_\e\r\|^2_{\mathcal{V}^D}\leq max\{C,1\}\l\|\psi_\e\r\|^2_{\mathcal{V}^{U^{\e}}}+C.$
In fact,
\begin{equation*}
\begin{split}
\l\|\psi_\e\r\|^2_{\mathcal{V}^D}&=\l\|\psi_\e\r\|^2_{L^2(D)}+C\int_{D\times \mathbb{R}^d}(\mathcal{D}^*\psi_\e)^2dxdz
=\l\|\psi_\e\r\|^2_{L^2(U^{\e})}+C\int_{U^{\e}\times \mathbb{R}^d}(\mathcal{D}^*\psi_\e)^2dxdz\\
&+C\int_{D/U^{\e}}\int_{\l(D/U^{\e}\r)^c}(m^{\e}(x))^2K^2(x,\frac{x}{\e})\nu(x,z)dzdx\\
&\leq max\{C,1\}\l\|\psi_\e\r\|^2_{\mathcal{V}^{U^{\e}}}+C \int_{D/U^{\e}}\int_{\l(D/U^{\e}\r)^c}(m^{\e}(x))^2K^2(x,\frac{x}{\e})\nu(x,z)dzdx.
\end{split}
\end{equation*}
From the fact that the derivative of the function $m^{\e}$ at $\{x\in D:dist (x,\partial D)=2\varepsilon\}$ is $0,$
we have
\begin{equation}\label{yao2}
\l\|\psi_\e\r\|^2_{\mathcal{V}^D}\leq max\{C,1\}\l\|\psi_\e\r\|^2_{\mathcal{V}^{U^{\e}}}+C.
\end{equation}
Finally, we will get the convergence rate.
 Combining (\ref{yao3}), (\ref{yao1}), (\ref{yao2}) and \cite[Corollary 5.1]{RK} we conclude
\begin{equation*}
\begin{split}
\l\|K_\e\r\|^2_{\mathcal{X}^D}&=\l\|\psi_\e\r\|^2_{\mathcal{X}^D}+C\leq C\l\|\psi_\e\r\|^2_{\mathcal{V}^D}+C
\leq  C\l\|\psi_\e\r\|^2_{\mathcal{V}^{U^\e}}+C
\leq C \varepsilon^{-\a}.
\end{split}
\end{equation*}
We thus estimate from Lemma \ref{le}
\begin{equation*}
\begin{split}
\l\|Z_\e\r\|_{\mathcal{V}^D}&\leq C\varepsilon^{\frac{1+\a}{2}}\l\|F_\e\r\|_{H^{-\a/2}(D)}+C\varepsilon^{\frac{1+\a}{2}}\l\|K_\e\r\|_{\mathcal{X}^D}
\leq C\varepsilon^{\frac{1+\a}{2}}+\varepsilon^{\frac{1+\a}{2}}\varepsilon^{-\frac{\a}{2}}c_{11}\leq C\varepsilon^{\frac{1}{2}}.
\end{split}
\end{equation*}
This completes the proof of Theorem \ref{T1}.


\begin{thebibliography}{119}
\bibitem{Duan}
Duan, J. and Wang, W.
Effective dynamics of stochastic partial differential equations.
Elsevier (2014).
\bibitem{Duant}
Duan, J. An introduction to stochastic dynamics, Cambridge University Press. Elsevier (2015).

\bibitem{DS}
Dipierro, S., Ros-Oton, X. and Valdinoci, E. Nonlocal problems with Neumann boundary conditions. Revista Matematica Iberoamericana 33(2), 377-416 (2017).
\bibitem{FM}
Felsinger, M., Kassmann, M. and Voigt, P. The Dirichlet problem for nonlocal operators. Mathematische Zeitschrift, 279(3-4), 779-809 (2015).
\bibitem{MV}
Millot, V., Sire, Y. and Wang, K. Asymptotics for the fractional Allen-Cahn equation and stationary nonlocal minimal surfaces. Archive for Rational Mechanics and Analysis, 231(2), 1129-1216 (2019).
\bibitem{RA}
Rutkowski, A. The dirichlet problem for nonlocal l\'{e}vy-type operators. Publicacions Matematiques, 62(1)(2017).
\bibitem{RK}
Bogdan, K., Grzywny, T., Pietruska-Paluba, K. and Rutkowski, A.  Extension and trace for nonlocal operators. Journal de Math¨¦matiques Pures et Appliqu¨¦es (2019).
%\bibitem{FM2}
%Felsinger, M., Kassmann, M. and Voigt, P. (2015). The Dirichlet problem for nonlocal operators. Mathematische Zeitschrift, 279(3-4), 779-809.
\bibitem{MW}
McLean, W. Strongly elliptic systems and boundary integral equations. Cambridge University Press (2000).

\bibitem{BK}
Bogdan, K., Dyda, B. and Luks, T.  On hardy spaces of local and nonlocal operators. Hiroshima Mathematical Journal, 44(2), 193-215 (2011).



\bibitem{Oleinik}
Oleinik, O. A., Shaniaev, A. S. and Yosifian, G. A.  Mathematical Problems in Elasticity and Homogenization. North-Holland, Amsterdam (1992).

\bibitem{GT}
Grzywny, T., Kassmann, M. and Lezaj, L.  Remarks on the nonlocal Dirichlet problem. Potential Analysis (2020).
\bibitem{11}
Chen, H., Duan, J. and  Lv, G. Boundary Blow-up Solutions to Nonlocal Elliptic Systems of Cooperative Type. Annales Henri Poincar\'{e} 19, 2115-2136 (2018).
\bibitem{12}
Chen, G. and Wang, Y.
Invariant measure of stochastic fractional Burgers equation with degenerate noise on a bounded interval. Communications on Pure and Applied Analysis 18(6):3121-3135(2019).

\bibitem{13}
Blomker, D., Hairer, M. and Pavliotis, G. A. Multiscale analysis for stochastic partial differential equations with quadratic nonlinearities. Nonlinearity, 20(7):1721-1744(2007)

\bibitem{GP}
Pavliotis, G. A. and Stuart, A. M. Multiscale Methods: Averaging and Homogenization. Springer. (2007)

\bibitem{DQ}
Du, Q., Gunzburger, M. and Lehoucq, R. B., et al. Analysis and approximation of nonlocal diffusion problems with volume constraints. SIAM Review, 54(4):667-696(2012).
\end{thebibliography}
\end{document}